\documentclass[11pt]{article}
\usepackage[a4paper,width=6 in, height=8.8 in]{geometry}
\usepackage{amsmath,amssymb}
\usepackage{hyperref}
\usepackage{setspace}
\usepackage[utf8]{inputenc}
\usepackage[english]{babel}
\usepackage{booktabs,multirow}
\usepackage{blkarray}
\usepackage{lscape}
\usepackage[T1]{fontenc}
\usepackage{authblk,color}
\usepackage[english]{babel}
\usepackage{amsfonts, amstext, amsthm, booktabs, dcolumn}
\usepackage{graphicx} 
\usepackage{float}
\usepackage{caption}
\usepackage{subcaption}
\usepackage[mathscr]{euscript}
\usepackage{pgfplots}
\pgfplotsset{compat=1.7}

\newtheorem{theorem}{\bf Theorem}[section]

\newtheorem{corollary}{\bf Corollary}[section]
\newtheorem{remark}{\bf Remark}[section]
\newtheorem{remarks}{\bf Remarks}[section]

\newtheorem{example}{\bf Example}[section]

\newcommand\floor[1]{\left\lfloor#1\right\rfloor}

\numberwithin{equation}{section}
\begin{document}
\title{Poisson Approximation to Weighted Sums of Locally Dependent Random Variables}
\author[ ]{Pratima Eknath Kadu}
\affil[ ]{\small Department of Maths \& Stats,}
\affil[ ]{\small K J Somaiya College of Arts and Commerce,}
\affil[ ]{\small Vidyavihar, Mumbai-400077.}
\affil[ ]{\small Email: pratima.kadu@somaiya.edu}
\date{}
\maketitle

\begin{abstract}
\noindent 
This paper deals with Poisson approximation to weighted sums of locally dependent random variables using Stein's method. The derived result represent a significant improvement of existing results. To illustrate the effectiveness of our findings, a numerical comparison is provided.
\end{abstract}

\noindent
\begin{keywords}
Poisson distribution; weighted sums; Stein's method.
\end{keywords}\\
{\bf MSC 2020 Subject Classifications:} Primary: 62E17, 60F05; Secondary: 62E20, 60E05.

\section{Introduction and Preliminaries}\label{sec:1}
Weighted sums of random variables (rvs) have been successfully applied in the areas related to least-square estimators, non-parametric regression function estimators, and Jackknife estimates, among many others. These sums also play a crucial role in characterizing various statistics. However, determining the exact distributions of weighted sums of rvs is generally impractical, prompting interest in approximations to study their asymptotic behaviour. The literature has extensively explored limit theorems for weighted sums of rvs, as evidenced by works such as  Chow and Lai \cite{CL1973}, Dolera \cite{DE2013}, Lai and Robbins \cite{LR1978}, Zhang \cite{ZL1989}, among many others.\\
In this paper, we consider the locally dependent set up similar to Kumar \cite{k2021}. Let $Y_1,Y_2,\ldots,$ $Y_n$ be $\mathbb{Z}_{+}$-valued random variables. For each $i$, let $c_i \in \mathbb{N}=\{1,2,\ldots\}$ and $i \in \mathcal{A}_i \subset \{1,2,\ldots,n\}$ be such that $Y_i$ is independent of $Y_{\mathcal{A}_i^c}$, where $Y_\mathcal{A}$ is the collection of random variables $\left\{Y_i, i \in \mathcal{A}\right\}$ and $\mathcal{A}^c$ denotes the complement of the set $\mathcal{A}$. Also, see Kumar and Vellaisamy \cite{KP2023}  and R\"{o}llin [19] for a similar type of locally dependent structure. In addition, if $\mathcal{A}_i=\{i\}$ then our locally dependent structure reduced to the independent collection of random variables. Define
\begin{align}
Z_n:= \sum_{i=1}^{n} c_i Y_{i,n}, \label{2:Zn}
\end{align}
the weighted sums of locally dependent rvs. We assume $c_i=1$ for at least one $i$ so that $Z_n$ becomes $\mathbb{Z}_+$-valued random variable (rv). Throughout the paper, let $Z_A:=\sum_{j\in A}c_iY_i$. Limit theorems for the weighted sums of rvs have been discussed in the literature if the sum of weights is finite. For example, Chow and Lai \cite{CL1973} are considered the finite sum of the square of weights, and Bhati and Rattihalli \cite{BR2014} are considered the geometrically weighted sums. But if the sums of weights are not finite then the study of limit behavior of such distributions becomes challenging. Therefore, in this paper, we consider natural weights and obtain approximations for $Z_n$ that helps to visualize the behavior of $Z_n$.\\
Next, let $X \sim \mathcal{P}_{\lambda}$, the Poisson distribution, then the probability mass function of $X$ is given by
\begin{align}
\mathbb{P}(X=k)=\frac{e^{-\lambda}{\lambda}^{k}}{k!}, \quad k\in \mathbb{Z}_+.\label{2:pmf}
\end{align}
We obtain error bound for X-approximation to $Z_n$. The total variance distance is used as a distance metric. Our results are derived using Stein's method (Stein \cite{stein1972}) which can be followed in three steps. First, for a random variable $Y$, find a Stein operator $\mathscr{A}$ such that $\mathbb{E}(\mathscr{A}g(Y))=0$, for $g \in G_Y$, where $G_Y:=\{g \in G|g(0)=0, ~\text{and}~ g(k)=0 ~\text{for}~ k \notin S(Y)\}$, $G:=\{f|f:\mathbb{Z}_{+}\rightarrow \mathbb{R} \text{ is bounded}\}$, and $S(Y)$ is the support of random variable $Y$. Second, solve the following Stein equation
\begin{equation}
\mathscr{A}g(k)=f(k)-\mathbb{E}f(Y),~ k \in \mathbb{Z}_+ ~ {\rm and} ~ f \in G. \label{2:StEq}
\end{equation}
Finally, putting a random variable $Z$ in place of $k$ and taking expectation and supremum, we have
\begin{equation}
d_{TV}(Y,Z):=\sup_{f \in \mathcal{I}} \big| \mathbb{E}f(Y)- \mathbb{E}f(Z) \big| =\sup_{f \in \mathcal{I}} \big| \mathbb{E}[\mathscr{A}g(Z)] \big|, \label{dist-xy}
\end{equation}
where $\mathcal{I}=\{{\bf 1}_A~|~A\subseteq \mathbb{Z}_+\}$ and ${\bf 1}_A$ is the indicator function of $A$. Now, for the random variable X defined in \eqref{2:pmf}, the Stein operator is given by
\begin{equation}
\mathscr{A}g(k)=\lambda g(k+1)-kg(k), \quad {\rm for}~ k \in \mathbb{Z}_+ ~{\rm and}~ g \in G_X \label{2:sop}
\end{equation}
and the solution to the Stein equation \eqref{2:StEq}, say $g_f$, satisfies
\begin{equation}
\|\Delta g_f\| \leq  \frac{1}{\max(1,\lambda)},\quad \text{for $f\in \mathcal{I}$, $g\in G_X$}. \label{2:bndX}
\end{equation}
For additional details, see Barbour and Hall \cite{BH}, Barbour {\em et. al.} \cite{BHJ}, Chen {\em et. al.} \cite{CGS}, Brown and Xia \cite{BX}, and Kumar {\em et. al.} \cite{KUV2,KPV2020}. For recent developments, see Ley {\em et. al.} \cite{LRS}, Upadhye and Barman \cite{NK}, Kumar \cite{k2021}, and references therein.\\
This paper is organized as follows. In Section \ref{2:MRs}, we obtain the error bound for Poisson approximation to $Z_n$. We show that our bound is better than the existing bound by numerical comparisons.  

\section{Bounds For Poisson Approximation} \label{2:MRs}
In this section, we derive an error bound for the Poisson approximation to weighted sums of locally dependent rvs. Additionally, we present the result for the independent case. Recall that $\mathcal{P}_{\lambda} $ follows Poisson distribution with parameter $\lambda > 0$ and $Z_n = \sum_{i=1}^{n} c_i Y_{i}$, where $ Y_{i}$, $i=1 ,2,\ldots, n $, are locally dependent rvs. The subsequent theorem establishes an upper bound for the total variation distance between  $Z_n$ and $\mathcal{P}_{\lambda} $.

\begin{theorem} \label{2:thm}
Let $Z_n$ and $\mathcal{P}_{\lambda} $ be defined as in \eqref{2:Zn} and \eqref{2:pmf}, respectively. Then
\begin{align*}
d_{TV}(Z_n,\mathcal{P}_{\lambda}) \le \frac{1}{\max(1,\lambda)}\sum_{i=1}^{n}\left\{\sum_{j\in \mathcal{A}_i}c_i^2[\mathbb{E}(Y_iY_j)+\mathbb{E}(Y_i)\mathbb{E}(Y_j)]-c_i\mathbb{E}(Y_i)\right\},
\end{align*}
where $\lambda=\mathbb{E}(Z_n)$.
\end{theorem}
\begin{proof}
Replacing $k$ by $Z_{n}$ in \eqref{2:sop} and taking expectation, we have
\begin{align}
\mathbb{E} [\mathscr{A} g(Z_n)]&={\lambda} \mathbb{E} [g(Z_n+1)]- \mathbb{E} [Z_ng(Z_n)] \nonumber \\ 
&=\sum_{i=1}^{n} c_i \mathbb{E} (Y_{i}) \mathbb{E} [g(Z_n+1)]- \sum_{i=1}^{n} c_i \mathbb{E} \left[ Y_{i} g(Z_n) \right]. \label{2:e1}
\end{align}
Next, let $Z_{i}^*=Z_n-\sum_{j\in \mathcal{A}_i}c_i Y_{i}=Z_n-Z_{\mathcal{A}_i}$ then $Y_{i}$ and $Z_{i}^*$ are independent rvs. Adding and subtracting $\sum_{i=1}^{n} c_i \mathbb{E} \left[ Y_{i} g(Z_{i}^*+1) \right]$ in \eqref{2:e1}, we get 
\begin{align*}
\mathbb{E} [\mathscr{A} g(Z_n)]&=\sum_{i=1}^{n} c_i \mathbb{E} (Y_{i}) \mathbb{E} [g(Z_n+1)-g(Z_i^*+1)]- \sum_{i=1}^{n} c_i \mathbb{E} \left[ Y_{i} (g(Z_n) -g(Z_i^*+1))\right]\\
&=\sum_{i=1}^{n} c_i \mathbb{E} (Y_{i}) \mathbb{E} \left[\sum_{j=1}^{Z_{\mathcal{A}_i}}\Delta g(Z_i^*+j)\right]- \sum_{i=1}^{n} c_i\mathbb{E} \left[Y_i\sum_{j=1}^{Z_{\mathcal{A}_i}-1}\Delta g(Z_i^*+j)\right]
\end{align*}
Note that if $Y_i\ge 1$, then $Z_{\mathcal{A}_i}\ge 1$. Therefore,
\begin{align*}
|\mathbb{E} [\mathscr{A} g(Z_n)]|\le \|\Delta g\|\left\{\sum_{i=1}^{n}c_i[\mathbb{E}(Y_i)\mathbb{E}(Z_{\mathcal{A}_i})+\mathbb{E}(Y_i(Z_{\mathcal{A}_i}-1))]\right\}.
\end{align*}
Using \eqref{2:bndX} the result follows.
\end{proof}

\begin{corollary}\label{2:cor1}
Let $Z_n$ be the weighted sums of independent rvs and $\mathcal{P}_{\lambda} $ be defined as in \eqref{2:pmf} satisfies $\lambda=\mathbb{E}(Z_n)$. Then
\begin{align*}
d_{TV}(Z_n,\mathcal{P}_{\lambda}) \le \frac{1}{\max(1,\lambda)} \sum_{i=1}^{n}\left\{c_i^2[\mathbb{E}(Y_i^2)+(\mathbb{E}(Y_i))^2]-c_i \mathbb{E}(Y_i)\right\}.
\end{align*}
\end{corollary}
\begin{proof}
The result follows by letting $\mathcal{A}_i=\{i\}$ in Theorem \ref{2:thm}.
\end{proof}

\begin{remark}
\begin{enumerate}
\item[(i)] If $Y_{i} \sim Ber(p_{i})$ and $c_i=1$, for $i=1,2,\ldots,n$, then it can be easily verified that
\begin{align}
d_{TV}(Z_n,\mathcal{P}_{\lambda}) &\le \frac{\sum_{i=1}^{n} {p}^{2}_{i}}{ \max \left(1,\lambda \right)},\label{2:l1}
\end{align}
which is the bound also obtained by Barbour and Hall \cite{BH}, and an improvement over the bounds given by Le Cam \cite{CAM} and Kerstan \cite{KJ}.

\item[(ii)] Note that the bound given in Corollary \ref{2:cor1} is an improvement over the bound given in Corollary 2.1 of Kumar \cite{k2021}. In particular, if $Y_{i} \sim Ber(p_{i})$ then, from Corollary \ref{2:cor1}, we have 
\begin{align}
d_{TV}(Z_n,\mathcal{P}_{\lambda}) &\le \frac{\sum_{i=1}^{n} (c_i^2{p}^{2}_{i}+c_i(c_i-1)p_i)}{ \max \left(1,\lambda \right)}, \label{2:re1}
\end{align}
where $\lambda= \sum_{i=1}^n c_i p_{i}$. Also, from Corollary 2.1 of Kumar \cite{k2021}, we have
\begin{align}
d_{TV}(Z_n,\text{PB}(N,p)) &\le \frac{\gamma}{\floor{N}pq} \sum_{i=1}^n c_i \left( \sum_{\ell=1}^{c_i-1} |q-\ell q_{i}|p_{i} + q c_i p_{i}^2 \right), \label{2:re2}
\end{align}
where PB$(N,p)$ follows pseudo-binomial distribution with $N=(1/p) \sum_{i=1}^n c_i p_{i}$ and $q=1-p=\left(\sum_{i=1}^n c_i^2 p_{i} q_{i}\right)/ \left(\sum_{i=1}^n c_i p_{i}\right)$. Also, $\gamma \le \sqrt{\frac{2}{\pi}} \big( \frac{1}{4} + \sum_{i=1}^n \gamma_i - \gamma^{*}\big)^{-1/2}$ with $\gamma_i= \min \{ 1/2, 1- d_{TV}(c_i Y_{i},c_i Y_{i}+1) \}$ and $\gamma^{*}= \max_{1 \le i \le n} \gamma_i$. Note that the bound given in \eqref{2:re1} is better than the bound given in \eqref{2:re2}. Also, the bound given in \eqref{2:re2} seems to be invalid ($q \ge 1$) for small values of $p_{i}$. For example, let the values of $p_{i}$ be as follows:
\begin{table}[H]
  \centering
  \caption{The values of $p_{i}$}
  \label{2:tab1}
  \begin{tabular}{|c|ccc|ccc|ccc|}
\hline
& $i$ & $p_{i}$ & $c_i$ & $i$ & $p_{i}$ & $c_i$ &  $i$ & $p_{i}$ & $c_i$ \\
\hline
Set 1&1-10 & 0.5 & 1&11-20 & 0.45 & 2 & 21-30 & 0.40 & 1\\
\hline
Set 2 & 1-10 & 0.05 & 1&11-20 & 0.04 & 2 & 21-30 & 0.04 & 3\\
\hline
\end{tabular}
\end{table}
\noindent
Then, the following table gives the comparison between the bounds given in \eqref{2:re1} and \eqref{2:re2}.
\begin{table}[H]
  \centering
  \begin{tabular}{|l|ll|ll|lllll}
\hline
\multirow{2}{*}{$n$} & \multicolumn{2}{|c|}{Set 1} &  \multicolumn{2}{|c|}{Set 2} \\
\cline{2-5}
 & From \eqref{2:re1} & From \eqref{2:re2}  & From \eqref{2:re1} & From \eqref{2:re2}\\
\hline
10  & 0.5 & 0.797885  & 0.025 & 0.797885\\
20  & 1.4 & 1.60360   & 0.683846 & Not Valid\\
30  & 1.17778 & 1.34907  & 1.3732 & Not Valid\\
\hline
\end{tabular}
\end{table}
\noindent
Observe that our bounds are better than the bounds given by Kumar \cite{k2021}.
\end{enumerate}
\end{remark}

\bibliographystyle{PV}
\bibliography{PEK}

\end{document}